\documentclass[11pt,leqno]{amsart}
\usepackage{a4wide}
\usepackage{amssymb,upgreek}
\usepackage{amsmath}
\usepackage{amsthm}
\usepackage{mathrsfs}
\usepackage[utf8]{inputenc}
\usepackage{bm,euscript,csquotes,comment,combelow}
\usepackage{tikz-cd}

\usepackage[backend=biber,style=alphabetic]{biblatex}
\usepackage[dvipsnames]{xcolor}
\usepackage{hyperref}
\usepackage{bbm}
\theoremstyle{plain}

\makeatletter
\def\cal@symb#1|#2{\expandafter\def\csname #2#1\endcsname{\mathcal{#1}}}
\def\calsymbols#1#2{\@for\@tmpz:=#2\do{\expandafter\cal@symb\@tmpz|{#1}}}
\def\bb@symb#1|#2{\expandafter\def\csname #2#1\endcsname{\mathbb{#1}}}
\def\bbsymbols#1#2{\@for\@tmpz:=#2\do{\expandafter\bb@symb\@tmpz|{#1}}}
\def\bold@symb#1|#2{\expandafter\def\csname #2#1\endcsname{\mathbf{#1}}}
\def\boldsymbols#1#2{\@for\@tmpz:=#2\do{\expandafter\bold@symb\@tmpz|{#1}}}
\def\scr@symb#1|#2{\expandafter\def\csname #2#1\endcsname{\mathscr{#1}}}
\def\scrsymbols#1#2{\@for\@tmpz:=#2\do{\expandafter\scr@symb\@tmpz|{#1}}}
\def\frak@symb#1|#2{\expandafter\def\csname #2#1\endcsname{\mathfrak{#1}}}
\def\fraksymbols#1#2{\@for\@tmpz:=#2\do{\expandafter\frak@symb\@tmpz|{#1}}}

\def\dmth@p#1|{\expandafter\let\csname#1\endcsname\relax
  \expandafter\DeclareMathOperator\csname#1\endcsname{#1}}
\def\operators#1{\@for\@tmpz:=#1\do{\expandafter\dmth@p\@tmpz|}}
\makeatother
\calsymbols{c}{A,B,C,D,E,F,G,H,I,J,K,L,M,N,O,P,Q,R,S,T,U,V,W,X,Y,Z}
\bbsymbols{b}{A,B,C,D,E,F,G,H,I,J,K,L,M,N,O,P,Q,R,S,T,U,V,W,X,Y,Z}
\boldsymbols{bf}{A,B,C,D,E,F,G,H,I,J,K,L,M,N,O,P,Q,R,S,T,U,V,W,X,Y,Z}
\scrsymbols{s}{A,B,C,D,E,F,G,H,I,J,K,L,M,N,O,P,Q,R,S,T,U,V,W,X,Y,Z}
\fraksymbols{fr}{a,b,c,d,e,f,g,h,i,j,k,l,m,n,o,p,q,r,s,t,u,v,w,x,y,z,X,B}

\newcommand{\zagrade}[1]{\left[ {#1}\right]}

\newcommand{\Cp}{\mathbb{C}_p}

\newcommand{\B}{\mathbb{B}}
\newcommand{\ol}{o_L}

\newcommand{\lt}{\operatorname{LT}}
\newcommand{\ox}{\mathscr{O}\left( \x \right)}
\newcommand{\oxb}{\mathscr{O}^b\left( \x \right)}
\newcommand{\ob}{\mathscr{O}\left( \B \right)}

\renewcommand{\sup}{\operatorname{sup}}
\newcommand{\xpm}{\x\zagrade{\pi^{>M}}}
\newcommand{\cyc}{\chi_{\text{cyc}}}
\newcommand{\tr}{\operatorname{Tr}}
\newcommand{\gal}{\operatorname{Gal}}
\renewcommand{\int}{\operatorname{Int}}
\renewcommand{\d}{\operatorname{D}\left(\ol, L\right)}
\newcommand{\db}{\operatorname{D}^b\left( \ol, L\right)}
\newcommand{\x}{\mathfrak{X}}
\newcommand{\xpinf}{\mathfrak{X}\zagrade{p^{\infty}}}
\newcommand{\la}{\lambda}
\newcommand{\xb}{\Lambda_L\left(\x\right)}
\newcommand{\clan}{\operatorname{C}^{L-\text{an}}\left(\ol, \Cp\right)}

\renewcommand{\hom}{\operatorname{Hom}}
\newcommand{\xpn}{\x\zagrade{\pi^{>N}}}
\newcommand{\lolpm}{\left(L/\ol\right)\zagrade{\pi^{>M}}}

\theoremstyle{definition}
\newtheorem{definition}{Definition}[section]
\newtheorem{theorem}[definition]{Theorem}
\newtheorem{corollary}[definition]{Corollary}
\newtheorem{proposition}[definition]{Proposition}

\newtheorem{lemma}[definition]{Lemma}

\newtheorem{remark}[definition]{Remark}

\usepackage{blkarray}
\usepackage{comment}
\usepackage{enumerate}
\usepackage{float}
\usepackage{graphicx}
\usepackage{listings}

\title{On the ring of bounded functions on the character variety}
\author{Uro\v{s} A.\ Colovi\'{c}}

\address{Uro\v{s} A.\ Colovi\'{c}, Mathematical Institute, University of Oxford}
\email{uros.colovic@maths.ox.ac.uk}

\date{12 November 2026}

\begin{document}

\begin{abstract}
Let $\mathfrak{X}$ denote the character variety associated to the group $(\ol, +)$, where $\ol$ is the ring of integers in a finite extension $L/\mathbb{Q}_p$. $\mathfrak{X}$ is a rigid analytic variety parametrizing locally $L$-analytic characters. Let $\xb$ denote the global sections bounded by $1$. In this paper we prove that $\xb = \ol[[\ol]]$, answering the question posed in \cite{bounded}.
\end{abstract}


\maketitle

\tableofcontents

\section{Introduction}

Let $L$ be a finite extension of $\mathbb{Q}_p$ with completed algebraic closure $\mathbb{C}_p$. The character variety $\mathfrak{X}$ was first defined by Schneider and Teitelbaum \cite{padic_fourier_theory} in the context of $p$-adic Fourier theory. $\mathfrak{X}$ is a rigid analytic variety defined over $L$ whose set of $K$-points $\mathfrak{X}\left(K\right)$, for $L \le K\le \Cp$ complete, is the set of locally $L$-analytic characters $\chi:(\ol, +)\to(K^{\times}, \cdot)$.

The character variety appears in the approach of Berger, Schneider and Xie to the Lubin-Tate generalization of cyclotomic $(\varphi, \Gamma)$-modules in \cite{bergerschxie}. This approach is motivated by the possibility of generalizing the $p$-adic local Langlands correspondence from $\operatorname{GL}_2(\mathbb{Q}_p)$ to $\operatorname{GL}_2(L)$. One of the reasons that $\x$ appears is that the ring of global bounded functions, denoted $\oxb$, has a natural action of the monoid $\ol-\{0\}$ analogous to that in the classical setting. Lubin-Tate analogues of cyclotomic $(\varphi, \Gamma)$-modules are defined over the rings that we get by localizing and completing $\oxb$.

In this paper we consider the subring $\xb \le \oxb$ consisting of global functions bounded by $1$ on the character variety. Ardakov and Berger in \cite{bounded} study, among other things, the question whether every  distribution bounded by $1$ comes from an integral measure, or equivalently, whether $\xb = \ol[[\ol]]$, see \cite[Question 1.1.1]{bounded}. Here $\ol[[\ol]]$ denotes the completed group ring and the map $\ol[[\ol]]\longrightarrow \xb$ is given by integration against the measure.

In \cite{bounded} several different criteria that imply $\xb = \ol[[\ol]]$ are proved. The question was further studied in the quadratic unramified case by Berger and Sprang in \cite{berger_sprang}.
In this paper we prove 
\begin{theorem}\label{cilj}
    It holds that $\xb = \ol[[\ol]]$.
\end{theorem}

In independent work, Berger, Sprang and Steingart prove the same result for quadratic extensions of $\mathbb{Q}_p$ using different methods.

The starting point of the proof of Theorem \ref{cilj} is a reduction to the integrality of a trace operator by using \cite[Theorem 3.1.9]{bounded}, see Lemma \ref{redukcija na psii}. We prove the integrality of this operator using a version of analytic boundary results from \cite{berger_analiticka_granica} in combination with a modification of Tate's classical result on absolute values of traces in deeply ramified field extensions given in Lemma \ref{procena za trag}. A key step in the proof is the connection between partial trace operators, which are used to factorize the trace operator from \cite[Theorem 3.1.9]{bounded}, and a normalized arithmetic trace, see Proposition \ref{parcijalni psi = trag}. For more details, see the sketch of the proof of Theorem \ref{cilj} in Section \ref{skica dokaza}.

Theorem \ref{cilj} has the following interesting corollary for integer valued polynomials, see \cite[Theorem 4.2.21]{bounded}.
Let $\int(\ol, \ol)$ denote the ring of integer valued polynomials on $\ol$. Denote the formal Lubin-Tate action\footnote{For some choice of Lubin-Tate formal group} by $[a]X$.
Recall from \cite[Section 1.5]{bounded} the definition $\operatorname{{Pol}}\left(\ol[[X]]\right) = \sum_{i,j\ge 0}\ol\sigma_{i,j}$ where $([a]X)^i = \sum_j \sigma_{i,j}(a)X^j$. Here $ \sigma_{i,j} \in \int(\ol, \ol)$.
\begin{corollary}\label{int = pol}
     The following equality holds: $\operatorname{{Pol}}\left(\ol[[X]]\right) = \int(\ol, \ol)$.
\end{corollary}

Related results were studied in \cite{rosendeshalit, prvitvrdio}, where it was proven that the space of coefficients $\left\{\sigma_{i,j}\right\}_{i,j\ge 0}$ generates $\int(\ol, \ol)$ as an algebra over $\ol$. Theorem B from \cite{berger_sprang} implies that  $\operatorname{{Pol}}\left(\ol[[X]]\right)$ is $p$-adically dense in $\int(\ol, \ol)$. Corollary \ref{int = pol} strengthens these results by showing that  $\left\{\sigma_{i,j}\right\}_{i,j\ge 0}$ generate all integer valued polynomials as an $\ol$-module.

Note that when $L = \mathbb{Q}_p$, the Lubin-Tate formal group can be taken to be the formal multiplicative group $\widehat{\mathbb{G}}_m$. In this case Corollary \ref{int = pol} recovers the classical fact that all integer valued polynomials over $\mathbb{Q}_p$ are $\mathbb{Z}_p$-linear combinations of binomial coefficients $\left\{\binom{x}{n}\right\}_{n\ge 0}$. In this regard, Corollary \ref{int = pol} can be interpreted as a statement that the “generalized” binomial coefficients span all integer valued polynomials. 

\begin{remark}
     In \cite[Remark 2.5]{colmez}, Colmez showed that a locally $L$-analytic $\Cp$-valued distribution on $\ol$ with bounded Fourier transform has order at most $d-1$, where $[L : \mathbb{Q}_p] = d$. If such a distribution $\la$ is defined over $L$, Theorem \ref{cilj} implies that $\la$ has order $0$.
\end{remark}

\subsection{Structure of the paper}
Section $2$ collects the necessary background and provides a sketch of the proof of Theorem \ref{cilj}. Section $3$ contains the proof of Theorem \ref{cilj}. Section \ref{appa} contains the proof of a result about unions of fibers being an analytic boundary. Section \ref{appb} contains the proof of our reformulation of a classical result of Tate about the behaviour of normalised traces in deeply ramified extensions. 

\subsection{Acknowledgements}
I want to thank Konstantin Ardakov for many valuable discussions, interesting insights and invaluable guidance and support.  I want to thank Laurent Berger and Andrew Graham for the discussion related to the question. I want to thank Konstantin Ardakov, Laurent Berger, James Newton, Branislav Prvulovi\'{c} and Vuka\v{s}in Stojisavljevi\'{c} for reading the first draft of the paper. 

\subsection{AI Declaration}
AI assistance consisted of suggesting the relation with the arithmetic trace and pointing to the operators $\psi_i$. No AI was used during the writing of the paper, but it was used to check the proofs and review the literature.

\section{Notation and background}

Fix a prime $p$. Let $\Cp$ denote the completed algebraic closure of $\mathbb{Q}_p$. Let the norm on $\Cp$ be normalized such that $|p| = p^{-1}$. Let $L/\mathbb{Q}_p$ be a finite extension of degree $d$ with ramification index $e$ and a uniformizer $\pi$. Denote by $q$ the cardinality of the residue field. Let $\{\mu_{p^n}\}_{n=1}^{\infty}$ be a system of primitive $p^n$-th roots of unity and denote  $\mathbb{Q}_{p,n} = \mathbb{Q}_p(\mu_{p^n})$, $L_n = L(\mu_{p^n})$, $L_{\infty} = \varinjlim_n L_n$.

\subsection{The character variety}
The character variety $\mathfrak{X}$ is a rigid analytic variety defined over $L$ whose set of $K$-points $\mathfrak{X}\left(K\right)$, for $L \le K\le \Cp$ complete, is the set of locally $L$-analytic characters $\chi:(\ol, +)\to(K^{\times}, \cdot)$. Denote the global sections on $\x$ by $\ox$, the bounded global sections by $\oxb$, and the power-bounded global sections by $\xb$. For the precise definition of the character variety see \cite[Section 2]{padic_fourier_theory}.

The set $\x\left(\Cp\right)$ is a group under pointwise multiplication of characters and we denote this multiplication by $\xi\cdot\eta$ for $\xi, \eta \in \x\left(\Cp\right)$. For $\xi \in \x\left(\Cp\right)$ we will use additive notation $k\xi$ to mean $\xi^k$. There is an action of $\ol$ on $\x\left(\Cp\right)$ defined by $a\zeta(x) = \zeta (ax)$. Let $a \in \ol$ and denote the $a$-torsion under this action by $\x\zagrade{a}$. Define $\xpinf = \varinjlim_m\x\zagrade{p^m}$ and let $\xpm$ denote the set of characters in $\xpinf$ whose $\pi$-order is bigger than $M$.

\begin{remark}\label{olizo}
    There exists a noncanonical isomorphism $\xpinf\cong L/\ol$ as $\ol$-modules. One can use \cite[Lemma 2.13]{bergerschxie}, or argue directly as follows. Take a family of characters $\{\xi_m\}_{m=0}^{\infty}$ such that $\pi\xi_{m+1}=\xi_m$ and $\xi_m$ has $\pi$-order equal to $m$, and consider the family of $\ol$-module isomorphisms $\theta_m:\ol/\pi^m\ol \to \x\zagrade{\pi^m}$ defined by sending $a \mapsto a\xi_m$. These maps are compatible, so we can take the inductive limit to get an $\ol$-module isomorphism $\theta = \varinjlim_m\theta_m: L/\ol \to \xpinf$. Note that we view $L/\ol = \varinjlim_m\ol/\pi^m$ where the transition maps are given by multiplication by $\pi$.
\end{remark}

\subsection{$p$-adic Fourier theory and Lubin-Tate isomorphism}\label{furijeova teorija}

Let $\d$ denote the set of  locally $L$-analytic distributions on the additive group $\ol$. These carry an action of $\ol-\{0\}$ defined by $(a_*\la)(f(x)) = \la(f(ax))$. One of the main results of \cite[Theorem 2.3]{padic_fourier_theory} is an isomorphism 
\begin{equation}\label{dist = global na char}
    \d \cong \ox.
\end{equation}
This enables us to work with distributions, and we will be using this isomorphism implicitly\footnote{For example, in Theorem \ref{stabilnost pri parcijalnim tragovima}}. Note that if $\la \in \d$, we can evaluate $\la$ at a locally $L$-analytic $\Cp$-valued function on $\ol$.

Let $\db$ denote the distributions corresponding to  $\mathscr{O}^b\left(\x\right)$ under this isomorphism. We call these bounded distributions. There is a multiplicative norm on $\db$ given by
\begin{equation*}
    |\la| = \underset{\xi \in \x\left(\Cp\right)}{\sup}|\la(\xi)|.
\end{equation*}
Distributions that have norm at most one under the isomorphism (\ref{dist = global na char}) correspond to $\xb$. For more details see \cite[p. 3 Remark, Corollary 1.23]{bergerschxie}.

Let $\lt$ denote a Lubin-Tate group associated to $\pi \in \ol$, and let $\B$ denote the open unit ball over $\Cp$. Denote the base change of the character variety to $\Cp$ by $\x_{\Cp}$. In \cite[Theorem 3.6]{padic_fourier_theory} it is proved that 
\begin{equation*}
    \B \cong \x_{\Cp}.
\end{equation*}

\subsection{Operators $\varphi$ and $\psi$}\label{definicija psi i fi}

We recall the operators $\varphi, \psi_q$ used in \cite{bounded}, and transport them to distributions via $p$-adic Fourier theory and the $\lt$ isomorphism. 
Define the operator $\varphi : \ob \longrightarrow \ob$ by $Z \to \zagrade{\pi}Z$, where $\zagrade{\pi}$ denotes the Lubin-Tate action. This is a finite locally free map of constant rank, so we can define the normalized trace $\psi_q: \ob \longrightarrow \ob$ as the unique map such that
\begin{equation*}
    \varphi\psi_q = \frac{1}{q}\tr_{\ob/\varphi\left(\ob\right)}.
\end{equation*}
From now on, by transport and descent, we view these two operators as acting on distributions. We call them $\varphi, \psi:\d\to\d$.
After transport and descent, we get that $\varphi$ acts on distributions as multiplication by $\pi$, or in other words $\varphi = \pi_*$. 
 The action of $\psi$ is given by 
\begin{equation*}
    \left(\psi\la\right)\left(f\right) = \la \left(M_f\right), 
\end{equation*}
where $M_f(x)$ is defined as $0$ on $\ol-\pi\ol$ and as $f\left(\frac{x}{\pi}\right)$ otherwise.
For the proof of these claims see \cite[Remark 4.1.4, Section 4.1.4]{schneidervenjakob}.

\subsection{Analytic boundary}

We slightly modify the definition of an analytic boundary from \cite{berger_analiticka_granica} by restricting the family of distributions to $L$-valued ones. We fix the norm on $\x\left(\Cp\right)$ by the $\lt$ isomorphism from Section \ref{furijeova teorija}.

\begin{definition}
    A subset $\{a_n\}_{n=0}^{\infty}$ of $\x (\Cp )$ that satisfies $|a_n|\to 1$ with $n\to \infty$ is called an analytic boundary if for all $\la \in \db$ we have
    \begin{equation*}
        |\la| = \underset{n}\sup \left|\la(a_n)\right|.
    \end{equation*}
\end{definition}

The following result will play an important role.
\begin{theorem}[Berger]\label{berger granica}
    Let $N$ be a nonnegative integer. Then $\xpn$ is an analytic boundary. In particular, for $\lambda \in \db$ we have $|\la| = \sup_{\xi \in \xpn}|\la(\xi)|$.
\end{theorem}
The proof follows by combining \cite[Theorem D]{berger_analiticka_granica} and \cite[Lemma 1.1.3]{berger_analiticka_granica}.
We record one useful lemma here.

\begin{lemma}\label{ne menja normu}
    Let $\lambda \in \d$ and $a \in \ol$ with $a\neq 0$. Then $\la \in \db$ if and only if $a_*\la \in \db$. Furthermore, if $\la \in \db$, then
    \begin{equation*}
        |\la| = |a_*\la|.
    \end{equation*}
    Specifically, if $|\varphi \la| \le 1$ then $|\la|\le 1$.
\end{lemma}

\begin{proof}
    The proof is immediate since multiplication by $a$ is surjective on $\x\left(\Cp\right)$ by \cite[Lemma 4.1.10]{schneidervenjakob}.
\end{proof}

\subsection{Partial trace operators}
Here we define partial trace operators that will be used to factorize $\psi$.
Since $\x\zagrade{p} \cong \ol/p\ol$ is a $\mathbb{Z}/p\mathbb{Z}$-module of rank $d$, there exist $d$ elements $\zeta_1,\ldots,\zeta_d$ such that they form a $\mathbb{Z}/p\mathbb{Z}$-basis of  $\x\zagrade{p}$. Let $\ker \zeta_i = K_i \le \ol$.

\begin{definition}
    For $1 \le i \le d$ define 
    \begin{equation*}
        \psi_i:\d \longrightarrow\d
    \end{equation*}
    by $\psi_i(\la) = \la(1_{K_i }\cdot)$.
\end{definition}

\subsection{Sketch of the proof}\label{skica dokaza}
  
Recall that $\{\mu_{p^n}\}_{n=1}^{\infty}$ is a system of primitive $p^n$-th roots of unity.
We want to prove that $\psi(\xb)\le \xb$, i.e. that the normalized trace operator is integral, see \cite[Theorem 3.1.9]{bounded}. The reduction to the integrality of the partial trace operators $\psi_i$ is formal, see Lemma \ref{redukcija na psii}. 

Fix from now on one index $i$ and a distribution $\la \in \xb$. To prove the integrality of the operator $\psi_i$, i.e. that $|\psi_i\la|\le 1$, we can use the theorem of Berger \cite[Theorem D]{berger_analiticka_granica} that says that the norm of $\psi_i\la$ can be seen as the supremum of $|(\psi_i\la)(\xi)|$ over all $\xi \in \xpinf$.

The main idea of the proof is that for sufficiently large integers $n$ and characters $\xi \in \xpinf$ satisfying $p^n\xi = 1, p^{n-1}\xi = \zeta_i$, we have that the value $\psi_i\la(\xi)$ can be seen as the normalized arithmetic trace for the totally ramified field extension $L(\mu_{p^n})/L(\mu_{p^{n-1}})$, see Proposition \ref{parcijalni psi = trag} (compare to \cite[p. 1023]{bergfarq}). There are good estimates for these traces due to Tate, see \cite[3.1, Corollary 2]{tejt}, which allow us to bound $|(\psi_i\la)(\xi)|$ giving us the integrality of the operators $\psi_i$.

Two issues remain. The first one is that extensions  $L(\mu_{p^n})/L(\mu_{p^{n-1}})$ for small $n$ need not satisfy the hypothesis of Tate's theorem (for example, these fields could be equal), so we need to be able to take $n$ large enough. The second issue is that for the trace formula to hold we need to restrict to characters $\xi \in \xpinf$ satisfying $p^{n-1}\xi = \zeta_i$. This is where Proposition \ref{analiticka granica veliko N} comes in, which is a version of \cite[Theorem D]{berger_analiticka_granica} that restricts $\xpinf$ to the family of characters such that $p^{n-1}\xi = \zeta_i$ for large enough $n$.

The geometric interpretation of the equality between the evaluation of $\psi_i\la$ at a torsion point and the normalized arithmetic trace is that the trace operator for finite étale (or finite locally free) maps commutes with base change. We apply this to the following unramified cover of degree $p$, where $\x_{\ker\zeta_i}$ is the character variety associated to the subgroup $\ker\zeta_i$ and the vertical arrow on the right is the restriction.

  \begin{equation*}
  \begin{tikzcd}[column sep = 2cm]
  \operatorname{Sp}L(\mu_{p^n})  \arrow{r}{\xi} \arrow[d] &  \x \arrow[d] \\
  \operatorname{Sp}L(\mu_{p^{n-1}})\arrow{r}{\xi|_{\ker{\zeta_i}}} & \x_{\ker\zeta_i}
  \end{tikzcd}
  \end{equation*}
To elaborate further, take a large enough integer $n$ and take a character $\xi \in \x\left(\Cp\right)$ of order $p^n$ such that $p^{n-1}\xi = \zeta_i$. Then $\xi$ has values in $L(\mu_{p^n})$, and $\xi|_{\ker\zeta_i}$ has values in $L(\mu_{p^{n-1}})$. The fiber of the image of $\xi$ in $\x_{\ker\zeta_i}\left(\Cp\right)$ is equal to $\{\xi \cdot \zeta_i^k\}_{k=0}^{p-1}$. But this set is exactly the set of Galois conjugates of $\xi$ under the action of $\gal\left(L(\mu_{p^n})/L(\mu_{p^{n-1}})\right)$. 
Now the partial trace operator $\psi_i$ is equal to averaging over the fiber, see Lemma \ref{drugi nacin}. However, the fiber consists of Galois conjugates as noted, so the averaging over the fiber is actually the normalized arithmetic trace.

The idea to connect arithmetic trace to evaluation at torsion points is already in the literature, see for example \cite[p. 1023]{bergfarq}.

\section{Proof of Theorem \ref{cilj}}

We start the proof of Theorem \ref{cilj} by reducing the statement to the integrality of the partial trace operators $\psi_i$.

\begin{lemma}\label{redukcija na psii}
    To prove Theorem \ref{cilj} it suffices to prove that $\psi_i\xb\le \xb$ for every $1\le i \le d$.
\end{lemma}
\begin{proof}
    We have 
    \begin{equation}\label{kompozicija parcijalnih tragova}
        \prod_{i=1}^{d}\psi_i = \varphi^e\psi^e,
    \end{equation}
    where the product means composition of operators. Indeed, we have that $\cap_{i=1}^dK_i = p\ol = \pi^{e}\ol$, and therefore 
    \begin{equation*}
        \left(\prod_{i=1}^{d} \psi_i \right)\la = \la(1_{p\ol}\cdot) = \varphi^e \psi^e (\la).
    \end{equation*}
    By  \cite[Theorem 3.1.9]{bounded}, to prove Theorem \ref{cilj} it suffices to prove $\psi(\xb)\le \xb.$ By the assumption we have  
    \begin{equation*}
        \left(\prod_{i=1}^{d}\psi_i\right)\xb\le \xb,
    \end{equation*}
    which, combined with (\ref{kompozicija parcijalnih tragova}), gives us that 
    \begin{equation*}
        \varphi^e\psi^e\xb\le\xb.
    \end{equation*}
    Lemma \ref{ne menja normu} gives us that 
    \begin{equation*}
        \psi^e\xb\le\xb.
    \end{equation*}
    Since $\varphi \xb\le\xb$ again by Lemma \ref{ne menja normu}, we get that 
    \begin{equation*}
        \psi^e\varphi^{e-1}\xb\le\xb.
    \end{equation*}
    Since $\psi\varphi = 1$, this finishes the proof.
\end{proof}

 Now we study partial trace operators further. First we give a different interpretation of $\psi_i$ as an average over the fiber under restriction $\x \to \x_{K_i}$. 

\begin{lemma}\label{drugi nacin}
    For $\la \in \d$ and $f \in \clan$ we have\footnote{By $\cdot$ we denote the pointwise multiplication of locally analytic functions, and recall that $k\zeta_i$ denotes $\zeta_i^k$}
    \begin{equation*}
        \psi_i\la(f) =\frac{1}{p} \sum_{k = 0}^{p-1}\la\left(f\cdot k\zeta_i\right).
    \end{equation*}
\end{lemma}

\begin{proof}
    The lemma is equivalent to the statement that
    \begin{equation*}
        \la(1_{K_i}\cdot f) =\frac{1}{p} \sum_{k = 0}^{p-1}\la\left(f\cdot k\zeta_i\right),
    \end{equation*}
    which follows because $\zeta_i$ is a character of order $p$ and this implies 
    \begin{equation*}
        1_{K_i} =\frac{1}{p} \sum_{k = 0}^{p-1} k\zeta_i.\qedhere
    \end{equation*}
    
\end{proof}

\begin{remark}
    In the Lubin-Tate coordinate, $\psi_i$ is given by $\psi_iF(Z) = \frac{1}{p}\sum_{k=0}^{p-1}F(Z\oplus \zagrade{k}\zeta_i)$.
\end{remark}

As a consequence of Lemma \ref{drugi nacin} we have that $\psi_i$ preserves bounded distributions
    \begin{equation*}
        \psi_i:\db \longrightarrow\db,
    \end{equation*}
and that these operators are bounded on $\db$. 

The idea for proving the integrality of the operators $\psi_i$ is testing their values against $\xi\in\xpinf$ by using Theorem \ref{berger granica}. Since such $\xi\in\xpinf$ are torsion, their values lie in cyclotomic extensions and hence they inherit a natural Galois action which we inspect. To bound the norm of $\psi_i\la$, for $\la \in \db$, we will interpret the value $\psi_i\la(\xi)$ as the value of the normalized trace of the cyclotomic extension. This equality is the key step in the entire proof. Recall that $\{\mu_{p^n}\}_{n=1}^{\infty}$ is a system of primitive $p^n$-th roots of unity and $L_n = L(\mu_{p^n})$,  $K_i = \ker \zeta_i \le \ol$.

\begin{proposition}\label{parcijalni psi = trag}
    Let $\xi \in \xpinf$ and let $n\ge 2$ be an integer. Suppose that $p^{n-1}\xi = \zeta_i$ and that $L_n/L_{n-1}$ is cyclic of degree $p$.  For $\la \in \d$ we have $\la(\xi) \in L_n$ and

    \begin{equation*}
        \psi_i\la(\xi) = \frac{1}{p}\tr_{L_n/L_{n-1}}\la(\xi).
    \end{equation*}
\end{proposition}

\begin{proof}
    Note that $p^n\xi = p\zeta_i = 1$, and hence $\xi$ takes values in $L_n$. Since torsion characters are locally constant, there exists a decomposition
    \begin{equation*}
        \xi = \sum_{\overline{a} \in \ol/p^n\ol}\xi(a)1_{a+p^n\ol}.
    \end{equation*}
    But because $\la$ is defined over $L$, this shows that $\la(\xi) \in L_n$.

    The group $G = \gal(L_n/L_{n-1})$ acts on $\x\zagrade{p^n}$, and we are interested in the orbit of $\xi$.
    First we prove that the orbit of $\xi$ has cardinality $p$. To do this, it is enough to prove that the stabilizer of $\xi$ has cardinality $1$. Therefore let $g \in G$ be such that $g\xi = \xi$. Note that since $\xi$ has exact order $p^n$, there exists an $a\in\ol$ such that $\xi(a)=\mu_{p^n}$. But we have that $g(\mu_{p^n}) = g\xi(a)=\xi(a) = \mu_{p^n}$, which implies that $g=1$.
    
    We claim that $\frac{g\xi}{\xi}\left(K_i \right) = 1$ for all $g \in G$. 
    To prove this, it suffices to prove that $G$ acts trivially on $\xi(K_i)$. From the assumption that $p^{n-1}\xi = \zeta_i$, we get that $\xi(K_i)^{p^{n-1}} = \left(p^{n-1}\xi\right)(K_i)=\zeta_i(K_i) = 1$. This implies that $\xi(K_i) \le L_{n-1}$, and so $G$ acts trivially on $\xi(K_i)$.
    
    We have proved that the set $\{\frac{g\xi}{\xi}\}_{g\in G}$ gives us $p$ different elements of $\hom(\ol/K_i, \Cp^{\times})$. On the other hand, since $\ol/K_i$ is cyclic of order $p$, we have that $\hom(\ol/K_i, \Cp^{\times}) = \{k\zeta_i\}_{k=0}^{p-1}$. This shows that $\left\{g\xi\right\}_{g\in G} = \left\{\xi \cdot k\zeta_i\right\}_{k=0}^{p-1}$.

    Note that $g\left(\la(\xi)\right) = \la(g\xi)$ because $\la$ is defined over $L$ and $\xi$ is torsion and hence locally constant. Therefore
    \begin{equation*}
        \frac{1}{p}\tr_{L_n/L_{n-1}}\la(\xi) = \frac{1}{p}\sum_{g \in \gal(L_n/L_{n-1})}g\left(\la(\xi)\right) = \frac{1}{p}\sum_{g \in \gal(L_n/L_{n-1})}\la(g\xi) = \frac{1}{p}\sum_{k=0}^{p-1}\la(\xi\cdot k\zeta_i).
    \end{equation*}
    Applying Lemma \ref{drugi nacin} finishes the proof.
\end{proof}

\begin{remark}
    Proposition \ref{parcijalni psi = trag} can also be deduced from the cyclotomic case of the trace formula in \cite[p. 1023]{bergfarq} by passing to the finite cyclic quotient defined by $\xi$ and extending scalars. The proof we give is direct.
\end{remark}

The next step is to bound this trace. It turns out this is a direct corollary of the classic result by Tate \cite[Section 3.1, Corollary 2]{tejt}. We prove it in Section \ref{appb}.

\begin{lemma}\label{procena za trag}
    For any $\varepsilon>0$ there exists an integer $N = N(\varepsilon)$ such that for every $n > N$ we have that $L_n/L_{n-1}$ is Galois cyclic totally ramified of degree $p$ and 
    \begin{equation*}
        \left|\frac{1}{p}\tr_{L_n/L_{n-1}}(x)\right|\le p^{\varepsilon},
    \end{equation*}
    for all $x \in o_{L_n}$.
\end{lemma}

Recall that we want to prove that $\left|\psi_i\la\right|\le 1$ for $\la \in \xb$.
At this point we would like to apply a description of the norm as a supremum of $|(\psi_i\la)(\xi)|$ over $\xi \in\xpn$, as given in Theorem \ref{berger granica}, and then apply Proposition \ref{parcijalni psi = trag} and Lemma \ref{procena za trag} to bound these values. However, in order to apply Proposition \ref{parcijalni psi = trag} and Lemma \ref{procena za trag}, we need to restrict to the family of $\xi \in \xpinf$ such that $p^{n-1}\xi=\zeta_i$ for large enough $n$. To this end, we prove the following strengthening of Theorem \ref{berger granica}.

\begin{proposition}\label{analiticka granica veliko N}
     Let $\zeta \in \xpinf$ and let $N$ be a positive integer. The set 
    \begin{equation*}
        B_{\zeta, N} = \{\xi \in \xpinf : \left( \exists m\in \mathbb{Z}, m> N \right),  p^m\xi = \zeta
        \}
    \end{equation*}
    is an analytic boundary. In particular, for $\lambda \in \db$ we have $|\la| = \sup_{\xi \in B_{\zeta, N}}|\la(\xi)|$.
\end{proposition}
The proof of this proposition uses Theorem \ref{berger granica} and is given in Section \ref{appa}.

We can now prove that the operators $\psi_i$ are integral by combining Proposition \ref{parcijalni psi = trag}, Lemma \ref{procena za trag} and Proposition \ref{analiticka granica veliko N}. 

\begin{theorem}\label{stabilnost pri parcijalnim tragovima}
   The operators $\psi_i:\d \rightarrow\d$ are integral. In other words, $\psi_i\xb\le \xb$.
\end{theorem}
\begin{proof}
    For $\la \in \xb$  it suffices to prove that $|\psi_i \la| \le p^{\varepsilon}$ for every $0<\varepsilon$. 
    Therefore, fix a nonzero $\la \in \xb$ and $\varepsilon>0$, and take $N$ to be large enough such that Lemma \ref{procena za trag} holds.  By Proposition \ref{analiticka granica veliko N} we have 
    \begin{equation*}
        |\psi_i\la| = \underset{\xi \in B_{\zeta_i, N}}{\sup}|\psi_i \la(\xi)|,
    \end{equation*}
    and therefore it is enough to prove that

    \begin{equation*}
        |\psi_i\la (\xi)|\le p^{\varepsilon}
    \end{equation*}
     for every $\xi \in B_{\zeta_i, N}$. In this regard, take $\xi \in B_{\zeta_i, N}$ and a positive integer $n$ such that we have $p^{n-1}\xi = \zeta_i$. By the definition of $B_{\zeta_i, N}$, we have that $n>N$.

    By Proposition \ref{parcijalni psi = trag}
    \begin{equation*}
        \psi_i\la(\xi) = \frac{1}{p}\tr_{L_n/L_{n-1}}\la(\xi),
    \end{equation*}
 and $\la(\xi) \in L_n$. Note also that since $|\la| \le 1$ it follows that $\la(\xi) \in o_{L_n}$. Now, by the way we chose $N$ to satisfy Lemma \ref{procena za trag}, and because $n>N$ we have 
 \begin{equation*}
     \left|\psi_i\la(\xi)\right| = \left|\frac{1}{p}\tr_{L_n/L_{n-1}}{\la(\xi)}\right|\le p^{\varepsilon}.
 \end{equation*}
 This finishes the proof.
 \end{proof}

\begin{remark}
    The fact that $L_n/L_{n-1}$ is totally ramified is important only for Tate's bound.
\end{remark}

The proof of Theorem \ref{cilj} follows by combining Theorem \ref{stabilnost pri parcijalnim tragovima} and Lemma \ref{redukcija na psii}.

\section{A union of fibers in $\xpinf$ is an analytic boundary}\label{appa}
In this section we prove Proposition \ref{analiticka granica veliko N}.
Fix a positive integer $N$ throughout this section. Let $x\in L/\ol$ be an element of $\pi$-order equal to $m$ and define
    \begin{equation*}
        B_{x, N} = \{y \in L/\ol : \left( \exists n\in \mathbb{Z}, n> N \right),  p^ny = x
        \}.
    \end{equation*}
    
Let $\lolpm$ be the set of elements of $L/\ol$ whose $\pi$-order is bigger than $M$.
The following lemma is just a linear algebra fact about the $\ol$-module $L/\ol$.
\begin{lemma}\label{pravljenjes}
    There exists a finite set $S = S_{x, N} \le \ol - \{0\}$ and an integer $M = M(m, N)$ such that for every $y \in \lolpm$ there exists $s \in S$ such that $sy \in B_{x, N}$.
\end{lemma}

\begin{proof}
     For $x = 0$ the proof is trivial, so we assume that $x\neq0$. Let $R\le\ol- \{0\}$ be a finite set such that every element of $\left(\ol/\pi^m\ol\right)^{\times}$ has a lift in $R$, and take 
    \begin{equation*}
        S = \{r\pi^i: r\in R, 0 \le i < e\}.
    \end{equation*}
    We can take $M$ large enough, and we will see that it suffices to take $M = 2eN + m$. We take $y \in \lolpm$, and note that, because of our choice of $M$, there exists a unique pair of integers $(l, i)$, with $0 \le i \le e-1$ and $l>N$, such that $p^l\pi^iy$ has $\pi$-order equal to $m$. Since the action of $\ol^{\times}$ on elements of $\pi$-order $m$ in $\ol/\pi^m$ is transitive, we can find an element $r \in R$ such that $rp^l\pi^iy=x$. Therefore $p^l(r\pi^i)y = x$, which implies that $r\pi^iy \in B_{x, N}$.
\end{proof}

\begin{remark}
    The set $S$ works for any $N$ and only depends on the $\pi$-order of $x$.
\end{remark}

\begin{proof}[Proof of Proposition \ref{analiticka granica veliko N}]
    For $\zeta = 1$ we have $B_{1,N}=\xpinf$, so by Theorem \ref{berger granica} we are done. Henceforth assume that $\zeta\neq1$. By Remark \ref{olizo} we have an isomorphism $ \xpinf\cong L/\ol$ as $\ol$-modules, so we can take $S$ and $M = M(m, N)$ as in Lemma \ref{pravljenjes}, where $m$ is the $\pi$-order of $\zeta$. By Theorem \ref{berger granica}, for any $\alpha \in \db$, we have
    \begin{equation*}
        |\alpha| = \sup_{\xi \in \xpm}|\alpha(\xi)|.
    \end{equation*}

    Let $0\neq \la \in \db$ and
    \begin{equation*}
        \la_S = \prod_{s\in S}(s_*\la).
    \end{equation*}
    Take some $\eta \in \xpm$. Then from Lemma \ref{pravljenjes} it follows that there exists some $s \in S$ such that $s\eta \in B_{\zeta, N}$. We have 
    \begin{equation*}
        |\la_S(\eta)| = |\la(s\eta)|
        \prod_{s' \in S-\{s\}}|\la(s'\eta)|
        \le \left(\underset{\xi \in B_{\zeta, N}}{\sup}|\la(\xi)| \right)
        |\la|^{|S|-1}.
    \end{equation*}
    Taking the supremum over $\eta \in \xpm$ and noting that $|\la_S| = |\la|^{|S|}$ by multiplicativity of the norm (Section \ref{furijeova teorija}) and Lemma \ref{ne menja normu}, we get
    \begin{equation*}
        |\la|^{|S|} = |\la_S| = \underset{\eta \in \xpm}{\sup}|\la_S(\eta)| \le |\la|^{|S|-1} \underset{\xi \in B_{\zeta, N}}{\sup}|\la(\xi)|.
    \end{equation*}
    Therefore we get 
    $|\la| \le \sup_{\xi \in B_{\zeta, N}}|\la(\xi)|$. The reverse inequality is obvious, showing that $|\la| = \sup_{\xi \in B_{\zeta, N}}|\la(\xi)|$. Finally, if $\xi \in B_{\zeta, N}$ has $\pi$-order $n$, by the proof of \cite[Theorem 3.5.4]{bounded} we have that $|\xi| = |\pi|^{\frac{1}{q^{n-1}(q-1)}}$, and the left-hand side tends to $1$ as $n$ tends to infinity, finishing the proof.
\end{proof}

\begin{remark}
    For the argument in Theorem \ref{stabilnost pri parcijalnim tragovima} we do not need that $B_{\zeta, N}$ has points whose norm is tending to $1$.
\end{remark}

\section{Normalized trace estimates for totally ramified $\mathbb{Z}_p$-extensions}\label{appb}

 In this section we prove Lemma \ref{procena za trag}. Recall the following result of Tate \cite[Section 3.1, Corollary 2]{tejt}. Here $K$ is a complete discretely valued field of characteristic $0$, with a perfect residue field of characteristic $p$.
\begin{theorem}[Tate]
    Let $K_{\infty}/K$ be a totally ramified Galois extension with Galois group $\mathbb{Z}_p$. Let $K_n$ correspond to the subgroup $p^n\mathbb{Z}_p$. Then for every integer $n\ge 1$, $K_n/K$ is cyclic of degree $p^n$ and there exists a constant $a$, independent of $n$, such that for any $x \in K_n$
    \begin{equation*}
        \left|\tr_{K_{n}/K_{n-1}}(x)\right|\le |x||p|^{1-ap^{-n+1}}.
    \end{equation*}
\end{theorem}

Recall that $\{\mu_{p^n}\}_{n=1}^{\infty}$ is a system of primitive $p^n$-th roots of unity and that  $\mathbb{Q}_{p,n} = \mathbb{Q}_p(\mu_{p^n})$, $L_n = L(\mu_{p^n})$, $L_{\infty} = \varinjlim_n L_n$.

\begin{proof}[Proof of Lemma \ref{procena za trag}]
     
     Denote the residue fields of $L_n$ by $k_n$, and note that the sequence $\zagrade{k_n:\mathbb{F}_p}$ is bounded and nondecreasing, and hence it stabilises. The sequence of numbers $\zagrade{L_n:\mathbb{Q}_{p,n}}$ is nonincreasing and hence it stabilises. 
      Take a large $m \ge 2$ such that both of the sequences stabilise. Then we have $\zagrade{L_n:L_{n-1}} = p$ for $n>m$ and these extensions are totally ramified.
     Next we prove that $L_{\infty}/L_n$ is totally ramified Galois with Galois group $\mathbb{Z}_p$ for all $n\ge m$.

     We have $\zagrade{L_{n+s}:L_n} = p^s$ and this extension is totally ramified. Since $L_{n+s} = L_n(\mu_{p^{n+s}})$, this is also Galois. Take $\sigma \in \gal\left(L_{n+s}/L_n\right)$. Then $\sigma $ fixes $\mu_{p^n}$ and is determined by what it does to $\mu_{p^{n+s}}$. This proves that the natural restriction map  
     \begin{equation*}
         \gal\left(L_{n+s}/L_n\right)\longrightarrow \gal\left(\mathbb{Q}_{p,n+s}/\mathbb{Q}_{p,n}\right)
     \end{equation*}
     is an injection. Since both groups have cardinality $p^s$, this has to be an isomorphism.
     Finally we get
     \begin{equation*}
         \gal\left(L_{\infty}/L_n\right)=\underset{s}{\varprojlim} \gal\left(L_{n+s}/L_n\right) \cong \underset{s}{\varprojlim}\gal\left(\mathbb{Q}_{p,n+s}/\mathbb{Q}_{p,n}\right) \overset{\cyc}{\underset{\sim}{\xrightarrow{\hspace*{1cm}}}} 1+p^n\mathbb{Z}_p,
     \end{equation*}
     and note that $1+p^n\mathbb{Z}_p \cong\mathbb{Z}_p$ since $n\ge 2$. This shows that $L_{\infty}/L_n$ is a totally ramified Galois extension with Galois group $\mathbb{Z}_p$ for all $n\ge m$.
     
     Finally let $K = L_m$ and $K_{\infty} = L_{\infty}$ and apply Tate's theorem to get
    \begin{equation*}
        \left|\frac{1}{p}\tr_{L_{m+n+1}/L_{m+n}}(x)\right| \le |x||p|^{-ap^{-n}} \le |p|^{-ap^{-n}} = p^{ap^{-n}}
    \end{equation*}
    for all $n$ and all $x \in o_{L_{m+n+1}}$. Now take $n$ large enough such that $p^{ap^{-n}}<p^{\varepsilon}$. Taking $N = n+m+1$ finishes the proof.
\end{proof}

\providecommand{\bysame}{\leavevmode\hbox to3em{\hrulefill}\thinspace}
\providecommand{\MR}{\relax\ifhmode\unskip\space\fi MR }
\providecommand{\MRhref}[2]{%
  \href{http://www.ams.org/mathscinet-getitem?mr=#1}{#2}
}

\printbibliography

\end{document}